\documentclass[12pt]{article}
\usepackage{amsfonts}
\usepackage{amsmath}
\usepackage{amssymb}
\usepackage{myart}

\normalbaselines
\font\tenbm=cmmib10
\font\sevenbm=cmmib7

\newfam\bmfam
\textfont\bmfam=\tenbm \scriptfont\bmfam=\sevenbm
\input tcilatex.tex
\begin{document}

\author{Yuri A. Rylov}
\title{Geometries with intransitive equivalence relation}
\date{Institute for Problems in Mechanics, Russian Academy of Sciences \\
101-1 ,Vernadskii Ave., Moscow, 117526, Russia \\
email: rylov@ipmnet.ru\\
Web site: {$http://rsfq1.physics.sunysb.edu/\symbol{126}rylov/yrylov.htm$}\\
or mirror Web site: {$http://gasdyn-ipm.ipmnet.ru/\symbol{126}%
rylov/yrylov.htm$}}
\maketitle

\begin{abstract}
One considers geometry with the intransitive equaivalence relation. Such a
geometry is a physical geometry, i.e. it is described completely by the
world function, which is a half of the squared distance function. The
physical geometry cannot be axiomatized, in general. It is obtained as a
result of deformation of the proper Euclidean geometry. Class of physical
geometries is more powerful, than the class of axiomatized geometries. The
physical geometry admits one to describe such geometric properties as
discreteness, granularity and limited divisibility. These properties are
important in application to the space-time. They admits one to explain the
discrimination properties of the space-time, which generate discrete
parameters of elementary particles. Mathematical formalism of a physical
geometry is very simple. The physical geometry is formulated in geometrical
terms (in terms of points and world function) without a use of means of
description (coordinate system, space dimension, manifold, etc.).
\end{abstract}

\section{Introduction}

Physical geometry is a science, which investigates mutual disposition of
geometrical objects in the space, or mutual disposition of events in the
event space (space-time). Geometrical objects are some subsets of points of
the point set $\Omega $, where the geometry is given. The point set $\Omega $
is called the space, or the space-time (if the points of $\Omega $ are
events). The physical geometry is the set of all propositions on properties
of all geometrical objects, belonging to the space $\Omega $. The physical
geometry $\mathcal{G}$ is determined completely, if the distance $\rho
\left( P,Q\right) $ between any pair of points $P,Q\in \Omega $ is given. It
is clear intuitively, that the mutual disposition of geometrical objects is
determined completely, if the distance $\rho \left( P,Q\right) $ is given
for all pairs of points. In the case of the proper Euclidean geometry the
fact, that the proper Euclidean geometry can be described completely in
terms of the distance $\rho $, has been proven \cite{R2001}. In general, the
world function $\sigma \left( P,Q\right) =\frac{1}{2}\rho ^{2}\left(
P,Q\right) $ is more convenient, because the world function is real in those
cases, when $\rho $ may be imaginary (for instance, in the geometry of
Minkowski). The term "world function" was introduced by J.L. Synge in his
presentation of the general relativity theory \cite{S60}.

The world function $\sigma $ is defined as follows%
\begin{equation}
\sigma :\qquad \Omega \times \Omega \rightarrow \mathbb{R},\qquad \sigma
\left( P,P\right) =0,\qquad \forall P\in \Omega  \label{b1.1}
\end{equation}

Construction of geomerical objects and description of their properties is
the main problem of any geometry. The set of geometrical objects, as well as
the set of all geometrical propositions, are continual sets, and it is very
difficult to construct them without a proper formalization.

The Euclidean method of the geometry construction admits one to construct
geometries, which can be axiomatized, i.e. such geometries, where the set $%
\mathcal{S}$ of all propositions can be deduced from the finite set $%
\mathcal{A}$ of basic suppositions (axioms) by means of the rules of the
formal logic. Axioms describe properties of basic geometrical objects
(point, segment, angle). Combining basic objects, one obtains more complex
geometrical objects, consisting of many blocks (basic objects).
Propositions, concerning these complex objects, are deduced by means of the
formal logic from the axioms, which describe the basic objects properties.

The proper Euclidean geometry $\mathcal{G}_{\mathrm{E}}$ can be constructed
by the Euclidean method. It has been shown, that the Euclidean axioms, used
for description of the basic objects (point, segment, angle) are consistent
\cite{H30}. Thus, the proper Euclidean geometry $\mathcal{G}_{\mathrm{E}}$
can be axiomatized. The possibility of axiomatization together with the
formal logic is used for multiplication of geometrical propositions, when
one obtains the continual set of all geometrical propositions from the
finite set of basic propositions (axioms).

As far as only one way of the geometrical propositions multiplication was
known, this way (deduction from axioms) was used for construction of
non-Euclidean geometries. One believes, that, changing in proper way the
geometry axiomatics (properties of basic objects), one can obtain another
(non-Euclidean) geometries. In order, that the construction of the
non-Euclidean geometry $\mathcal{G}$ by means of this method be possible,
two following conditions are to be fulfilled:

1. The axiomatics is consistent, i.e. any two different deductive methods of
the proposition deduction are to give the same result.

2. All propositions of the geometry $\mathcal{G}$ can be expressed in terms
of the world function $\sigma $ and only in terms of $\sigma $.

If the first condition is fulfilled, but the second one is not, the obtained
non-Euclidean geometry $\mathcal{G}$ is not a physical geometry. Such a
non-Euclidean geometry should be qualified as a mathematical geometry. Thus,
the mathematical geometry is such a set of propositions (logical
construction), which can be deduced from a system of axioms by means of the
rules of the formal logic. Some mathematical geometries do not satisfy the
second condition. Such geometries are not physical geometries. They cannot
be used for description of the space-time. For instance, the projective
geometry and the affine geometry are examples of a mathematical geometry,
which is not a physical geometry. The second condition is not fulfilled for
these geometries, because they do not contain the concept of distance (world
function).

However, the mathematical geometry may be a physical geometry
simultaneously. For instance, the proper Euclidean geometry $\mathcal{G}_{%
\mathrm{E}}$ is a mathematical geometry and a physical geometry
simultaneously, i.e. both conditions are fulfilled.

In a mathematical geometry a multiplication of the number of geometrical
propositions is produced by means of logical reasonings (proofs of theorems,
etc.). In a physical geometry, which is not a mathematical geometry
simultaneously, there are no means of multiplication of the number of
geometrical propositions, in general. In any physical geometry all
propositions are taken from some standard physical geometry, which is a
physical geometry and a mathematical geometry simultaneously. The set of
propositions in any physical geometry is the same as the set of propositions
in the standard physical geometry $\mathcal{G}_{\mathrm{st}}$. The standard
physical geometry $\mathcal{G}_{\mathrm{st}}$ is a physical geometry, which
can be axiomatized. Thus, the set of all geometrical propositions of $%
\mathcal{G}_{\mathrm{st}}$ is obtained by logical reasoning in the same way,
as in any mathematical geometry. Thereafter the set $\mathcal{S}$ of all
propositions $\mathcal{P}$ of $\mathcal{G}_{\mathrm{st}}$ is expressed in
terms of the world function $\sigma _{\mathrm{st}}$ of the geometry $%
\mathcal{G}_{\mathrm{st}}$ in the form $\mathcal{P=P}\left( \sigma _{\mathrm{%
st}}\right) $. It is possible, because $\mathcal{G}_{\mathrm{st}}$ is a
physical geometry.

The proper Euclidean geometry $\mathcal{G}_{\mathrm{E}}$ may be used as a
standard physical geometry $\mathcal{G}_{\mathrm{st}}$. Any proposition of $%
\mathcal{G}_{\mathrm{E}}$ may be expressed in terms of the Euclidean world
function $\sigma _{\mathrm{E}}$. For instance, let in the Cartesian
coordinate system of the $n$-dimensional proper Euclidean geometry the
vectors $\mathbf{P}_{0}\mathbf{P}_{1}$ and $\mathbf{Q}_{0}\mathbf{Q}_{1}$
have coordinates
\begin{equation}
\mathbf{P}_{0}\mathbf{P}_{1}=\left\{ x_{1},x_{2},...x_{n}\right\} ,\qquad
\mathbf{Q}_{0}\mathbf{Q}_{1}=\left\{ y_{1},y_{2},...y_{n}\right\}
\label{b1.2}
\end{equation}%
and the world function have the form
\begin{equation}
\sigma _{\mathrm{E}}\left( P_{0},P_{1}\right) =\frac{1}{2}%
\dsum\limits_{k=1}^{k=n}\left( x_{k}\right) ^{2}  \label{b1.3}
\end{equation}%
The scalar product $\left( \mathbf{P}_{0}\mathbf{P}_{1}.\mathbf{Q}_{0}%
\mathbf{Q}_{1}\right) $ of two vectors $\mathbf{P}_{0}\mathbf{P}_{1}$ and $%
\mathbf{Q}_{0}\mathbf{Q}_{1}$ is written in the form%
\begin{equation}
\left( \mathbf{P}_{0}\mathbf{P}_{1}.\mathbf{Q}_{0}\mathbf{Q}_{1}\right)
=\dsum\limits_{k=1}^{k=n}x_{k}y_{k}  \label{b1.4}
\end{equation}%
The scalar product (\ref{b1.4}) may be presented in terms of the Euclidean
world function $\sigma _{\mathrm{E}}$%
\begin{equation}
\left( \mathbf{P}_{0}\mathbf{P}_{1}.\mathbf{Q}_{0}\mathbf{Q}_{1}\right)
=\sigma _{\mathrm{E}}\left( P_{0},Q_{1}\right) +\sigma _{\mathrm{E}}\left(
P_{1},Q_{0}\right) -\sigma _{\mathrm{E}}\left( P_{0},Q_{0}\right) -\sigma _{%
\mathrm{E}}\left( P_{1},Q_{1}\right)  \label{b1.5}
\end{equation}%
Equivalency (equality) of vectors $\mathbf{P}_{0}\mathbf{P}_{1}$ and $%
\mathbf{Q}_{0}\mathbf{Q}_{1}$ is defined conventionally in the form
\begin{equation}
x_{k}=y_{k},\qquad k=1,2,...n  \label{b1.6}
\end{equation}

In terms of the world function the condition (\ref{b1.6}) of equivalency ($%
\mathbf{P}_{0}\mathbf{P}_{1}$eqv$\mathbf{Q}_{0}\mathbf{Q}_{1}$) of vectors $%
\mathbf{P}_{0}\mathbf{P}_{1}$ and $\mathbf{Q}_{0}\mathbf{Q}_{1}$ is written
in the form of two relations
\begin{eqnarray}
\mathbf{P}_{0}\mathbf{P}_{1}\mathrm{eqv}\mathbf{Q}_{0}\mathbf{Q}_{1}
&:&\qquad \left( \mathbf{P}_{0}\mathbf{P}_{1}.\mathbf{Q}_{0}\mathbf{Q}%
_{1}\right) =\left\vert \mathbf{P}_{0}\mathbf{P}_{1}\right\vert \cdot
\left\vert \mathbf{Q}_{0}\mathbf{Q}_{1}\right\vert  \label{b1.7} \\
\wedge \left\vert \mathbf{P}_{0}\mathbf{P}_{1}\right\vert &=&\left\vert
\mathbf{Q}_{0}\mathbf{Q}_{1}\right\vert  \label{b1.8}
\end{eqnarray}%
where
\begin{equation}
\left\vert \mathbf{P}_{0}\mathbf{P}_{1}\right\vert ^{2}=\left( \mathbf{P}_{0}%
\mathbf{P}_{1}.\mathbf{P}_{0}\mathbf{P}_{1}\right) =2\sigma _{%
\mathrm{E}}\left( P_{0},P_{1}\right)  \label{b1.9}
\end{equation}

It easy to verify that relation (\ref{b1.6}) is equivalent to relations (\ref%
{b1.7}), (\ref{b1.8}). It means that the conditions (\ref{b1.7}), (\ref{b1.8}%
) may be used as a definition of the equivalence relation of two vectors $%
\mathbf{P}_{0}\mathbf{P}_{1}$ and $\mathbf{Q}_{0}\mathbf{Q}_{1}$. According
to (\ref{b1.5}) the definition (\ref{b1.7}), (\ref{b1.8}) refers only to
points $P_{0},P_{1},Q_{0},Q_{1}$, defining vectors $\mathbf{P}_{0}\mathbf{P}%
_{1}$ and $\mathbf{Q}_{0}\mathbf{Q}_{1}$, and to the world function $\sigma
_{\mathrm{E}}$ between these points, whereas the definition (\ref{b1.6})
contains a reference to the dimension $n$ of the proper Euclidean geometry $%
\mathcal{G}_{\mathrm{E}}$ and to the Cartesian coordinate system in $%
\mathcal{G}_{\mathrm{E}}$. The definition (\ref{b1.6}) supposes a
possibility of introduction of a linear vector space in $\mathcal{G}_{%
\mathrm{E}}$. The condition (\ref{b1.7}) means that vectors $\mathbf{P}_{0}%
\mathbf{P}_{1}$ and $\mathbf{Q}_{0}\mathbf{Q}_{1}$ are parallel, whereas the
condition (\ref{b1.8}) means that their lengths are equal.

The definition (\ref{b1.7}), (\ref{b1.8}) is more general, than the
definition (\ref{b1.6}), because it does not use the means of description
(the coordinate system) and does not supposes introduction of the linear
vector space. The definition (\ref{b1.7}), (\ref{b1.8}) may be used in any
physical geometry $\mathcal{G}$, independently of whether or not one may
introduce the linear vector space in this geometry $\mathcal{G}$. Replacing
the Euclidean world function $\sigma _{\mathrm{E}}$ by the world function $%
\sigma $ of a physical geometry $\mathcal{G}$ in formulae (\ref{b1.6}) -- (%
\ref{b1.8}), one obtains the equivalence relation in the physical geometry $%
\mathcal{G}$. Such a replacement means a deformation of the proper Euclidean
geometry, when the distances between the points are changed.

The equivalence relation (\ref{b1.7}), (\ref{b1.8}) is transitive in the
proper Euclidean geometry, but it is intransitive, in general, in a physical
geometry. Indeed, if one fixes the points $P_{0},P_{1},Q_{0}$ and solves the
equations (\ref{b1.7}), (\ref{b1.8}) with respect to the point $Q_{1}$, one
obtains always one and only one solution in the proper Euclidean geometry.
This solution describes the well known property of the proper Euclidean
geometry, that at the point $Q_{0}$ there is one and only one vector $%
\mathbf{Q}_{0}\mathbf{Q}_{1}$, which is equal to vector $\mathbf{P}_{0}%
\mathbf{P}_{1}$ at the point $P_{0}$. Then for any points $%
P_{0},P_{1},Q_{0},Q_{1},R_{0}R_{1}$ it follows from $\mathbf{P}_{0}\mathbf{P}%
_{1}\mathrm{eqv}\mathbf{Q}_{0}\mathbf{Q}_{1}\wedge \mathbf{Q}_{0}\mathbf{Q}%
_{1}\mathrm{eqv}\mathbf{R}_{0}\mathbf{R}_{1}$, that $\mathbf{P}_{0}\mathbf{P}%
_{1}\mathrm{eqv}\mathbf{R}_{0}\mathbf{R}_{1}$. This property of the
equivalence relation is known as the transitivity. The transitivity of the
equivalence relation is a special property of the proper Euclidean geometry $%
\mathcal{G}_{\mathrm{E}}$, which may not exist in arbitrary physical
geometry.

In the general case of the physical geometry $\mathcal{G}$, when the
Euclidean world function $\sigma _{\mathrm{E}}$ in relations (\ref{b1.7}), (%
\ref{b1.8}) is replaced by another function $\sigma $, one can guarantee
neither existence of the vector $\mathbf{Q}_{0}\mathbf{Q}_{1}$ which is
equivalent to $\mathbf{P}_{0}\mathbf{P}_{1}$, nor its uniqueness.

If there are two different vectors $\mathbf{Q}_{0}\mathbf{Q}_{1},\mathbf{Q}%
_{0}\mathbf{Q}_{1}^{\prime },$ which satisfy the conditions $\mathbf{P}_{0}%
\mathbf{P}_{1}\mathrm{eqv}\mathbf{Q}_{0}\mathbf{Q}_{1}$ and $\mathbf{P}_{0}%
\mathbf{P}_{1}\mathrm{eqv}\mathbf{Q}_{0}\mathbf{Q}_{1}^{\prime }\wedge
\mathbf{Q}_{0}\mathbf{Q}_{1}\overline{\mathrm{eqv}}\mathbf{Q}_{0}\mathbf{Q}%
_{1}^{\prime }$, where $\overline{\mathrm{eqv}}$ means nonequivalence, then
the equivalence relation is intransitive. Intransitivity of the equivalence
relation is connected with the multivariance of the equivalence relation. By
definition, if at the point \ $Q_{0}$ there are several vectors $\mathbf{Q}%
_{0}\mathbf{Q}_{1},\mathbf{Q}_{0}\mathbf{Q}_{1}^{\prime },\mathbf{Q}_{0}%
\mathbf{Q}_{1}^{\prime \prime },...$, which are equivalent to the vector $%
\mathbf{P}_{0}\mathbf{P}_{1}$ at the point $P_{0}$, but they are not
equivalent between themselves, the equivalence relation is multivariant for
the point $Q_{0}$ and for the vector $\mathbf{P}_{0}\mathbf{P}_{1}$. The
equivalence relation may be zero-variant for the point $Q_{0}$ and for
vector $\mathbf{P}_{0}\mathbf{P}_{1},$ if at the point $Q_{0}$ there is no
vector, which is equivalent to vector $\mathbf{P}_{0}\mathbf{P}_{1}$. The
zero-variance may be considered as a special case of the multivariance,
although the zero-variance is connected with intransitivity in other way.

The case, when at the point $P_{0}$ there are two vectors $\mathbf{P}_{0}%
\mathbf{P}_{1}$ and $\mathbf{P}_{0}\mathbf{P}_{2}$, which are equivalent,
but they do not coincide%
\begin{equation}
\mathbf{P}_{0}\mathbf{P}_{1}\mathrm{eqv}\mathbf{P}_{0}\mathbf{P}_{2},\quad
P_{1}\neq P_{2}  \label{b1.10}
\end{equation}%
is also possible, in general. It is a special case of intransitive
(multivariant) equivalence. In other words, a coincidence of two vectors and
their equivalence are different properties, in general. Two coinciding
vectors are always equivalent, however, two equivalent vectors, having the
same origin, may not coincide.

The physical geometry is multivariant (and intransitive), in general. The
proper Euclidean geometry is an unique example of a single-variant physical
geometry. Even pseudo-Euclidean geometry is multivariant, in general.

As an example, let us consider the 4-dimensional geometry of Minkowski,
which is a pseudo-Euclidean geometry of index $1$. Let in the inertial
coordinate system the points $P_{0}=Q_{0},$ $P_{1}$, $Q_{1}$ have
coordinates
\begin{equation}
P_{0}=Q_{0}=\left\{ 0,0,0,0\right\} ,\qquad P_{1}=\left\{
x^{0},x^{1},x^{2},x^{3}\right\} ,\qquad Q_{1}=\left\{
y^{0},y^{1},y^{2},y^{3}\right\}  \label{b1.11}
\end{equation}%
In this coordinate system the world function has the form%
\begin{equation}
\sigma _{\mathrm{M}}\left( P_{1},Q_{1}\right) =\sigma _{\mathrm{M}}\left(
x,y\right) =\frac{1}{2}\left( \left( x^{0}-y^{0}\right) ^{2}-\left( \mathbf{x%
}-\mathbf{y}\right) ^{2}\right) =\frac{1}{2}x^{k}y_{k}  \label{b1.12}
\end{equation}%
Here $x^{k}$ and $x_{k}$ are contravariant and covariant coordinates
respectively of the point $P_{1}$. A summation from $0$ to $3$ is produced
over repeated contravariant and covariant indices. The vectors $\mathbf{P}%
_{0}\mathbf{P}_{1}$ and $\mathbf{Q}_{0}\mathbf{Q}_{1}$ have the form%
\begin{equation}
\mathbf{P}_{0}\mathbf{P}_{1}=\left\{ x^{0},\mathbf{x}\right\} =\left\{
x^{0},x^{1},x^{2},x^{3}\right\} ,\qquad \mathbf{Q}_{0}\mathbf{Q}_{1}=\left\{
y^{0},\mathbf{y}\right\} =\left\{ y^{0},y^{1},y^{2},y^{3}\right\}
\label{b1.12a}
\end{equation}%
Condition of equivalence (\ref{b1.7}), (\ref{b1.8}) of vectors $\mathbf{P}%
_{0}\mathbf{P}_{1}$ and $\mathbf{Q}_{0}\mathbf{Q}_{1}$ has the form%
\begin{eqnarray}
\frac{1}{2}y^{k}y_{k}+\frac{1}{2}x^{k}x_{k}-\frac{1}{2}\left(
x^{k}-y^{k}\right) \left( x_{k}-y_{k}\right) &=&x^{k}x_{k}  \label{b1.14} \\
x^{k}x_{k} &=&y^{k}y_{k}  \label{b1.15}
\end{eqnarray}%
Setting%
\begin{equation}
\alpha ^{k}=x^{k}-y^{k},\qquad k=0,1,2,3  \label{b1.16}
\end{equation}%
one obtains from (\ref{b1.14}), (\ref{b1.15})%
\begin{equation}
\alpha ^{k}\alpha _{k}=0,\qquad 2y^{k}\alpha _{k}=0  \label{b1.17}
\end{equation}%
If the vector $\mathbf{Q}_{0}\mathbf{Q}_{1}$ is timelike, i.e. $y^{k}y_{k}>0$%
, then there is the unique solution of equations (\ref{b1.17})
\begin{equation}
\alpha ^{k}=0,\qquad x^{k}=y^{k},\qquad k=0,1,2,3  \label{b1.18}
\end{equation}%
and this result coincides with the conventional definition of the vector
equivalence in the space-time of Minkowski
\begin{equation}
x^{k}=y^{k},\qquad k=0,1,2,3  \label{b1.18a}
\end{equation}

If the vector $\mathbf{Q}_{0}\mathbf{Q}_{1}$ is spacelike, i.e. $%
y^{k}y_{k}<0 $, the solution of (\ref{b1.17}) has the form%
\begin{equation}
x=\left\{ x^{0},\mathbf{x}\right\} ,\qquad x^{0}=y^{0}+\alpha ,\qquad
\mathbf{x=y+}\alpha \mathbf{n}  \label{b1.19}
\end{equation}%
where $\alpha $ is an arbitrary real number, $\mathbf{n}$ is the unit
3-vector, which forms the angle $\phi $ with the 3-vector $\mathbf{y}$. The
angle $\phi $ is determined by the relation
\begin{equation}
\cos \phi =\frac{y^{0}}{\left\vert \mathbf{y}\right\vert }  \label{b1.20}
\end{equation}%
The angle $\phi $ is real, if the vector $\mathbf{Q}_{0}\mathbf{Q}_{1}$ is
spacelike or null ($\left\vert y^{0}\right\vert \leq \left\vert \mathbf{y}%
\right\vert $). In particular, if the coordinate system is chosen in such a
way, that the 3-vector $\mathbf{y=}\left\{ y^{1},0,0\right\} $, one obtains%
\begin{equation}
\mathbf{Q}_{0}\mathbf{Q}_{1}=\left\{ y^{0},y^{1},0,0\right\} ,\qquad \mathbf{%
P}_{0}\mathbf{P}_{1}=\left\{ y^{0}+\alpha ,y^{1}+\alpha \frac{y^{0}}{y^{1}}%
,0,0\right\}  \label{b1.21}
\end{equation}%
where $\alpha $ is an arbitrary real number. We shall refer to the physical
geometry, described by the world function (\ref{b1.12}) as the $\sigma $%
-Minkowskian physical geometry. The term "geometry of Minkowski" remains for
the geometry constructed in the conventional way. Thus, the $\sigma $%
-Minkowskian geometry is multivariant with respect to spacelike and null
vectors, and it is single-variant with respect to timelike vectors. The
equivalence relation in the $\sigma $-Minkowskian geometry is intransitive.
It means, that the $\sigma $-Minkowskian geometry cannot be axiomatized,
because possibility of axiomatization of a geometry generates a transitive
equivalence relation and a single-variant definition of equivalence.

At the conventional approach, based on a use of the linear vector space, the
geometry of Minkowski is single-variant with respect to all vectors.
Single-variance is conditioned by a use of the equivalence relation (\ref%
{b1.18a}), which is formulated in terms of coordinates. A use of coordinates
may be replaced by four relations, written in terms of the world function.
For instance%
\begin{equation}
x^{k}=g^{kl}x_{l},\qquad x_{k}=\left( \mathbf{P}_{0}\mathbf{P}_{1}.\mathbf{OS%
}_{k}\right) ,\qquad k=0,1,2,3  \label{b1.22}
\end{equation}%
where metric tensor is defined via the world function (\ref{b1.12}) and the
relation (\ref{b1.5}) by means of relations
\begin{equation}
g_{kl}=\left( \mathbf{OS}_{k}.\mathbf{OS}_{l}\right) ,\qquad
g^{kj}g_{jl}=\delta _{l}^{k},\qquad k,l=0,1,2,3  \label{b1.23}
\end{equation}%
The four relations (\ref{b1.6}) take the form%
\begin{equation}
\left( \mathbf{P}_{0}\mathbf{P}_{1}.\mathbf{OS}_{k}\right) =\left( \mathbf{Q}%
_{0}\mathbf{Q}_{1}.\mathbf{OS}_{k}\right) ,\qquad k=0,1,2,3  \label{b1.24}
\end{equation}%
They contain a reference to five additional points $%
O,S_{0},S_{1},S_{2},S_{3} $, representing the coordinate system. Although
these relations are invariant with respect to a choice of five points $%
O,S_{0},S_{1},S_{2},S_{3}$, but they are sensitive to number of these
points, which depends on the dimension of the geometry of Minkowski.
Definition of two vectors equivalence in a physical geometry is to depend
only on points $P_{0},P_{1},Q_{0},Q_{1}$, which determine vectors $\mathbf{P}%
_{0}\mathbf{P}_{1},\mathbf{Q}_{0}\mathbf{Q}_{1}$ and on the world function
between these points. Dependence on additional points means that the
considered geometry is not a physical geometry. The considered geometry is a
fortified physical geometry, i.e. a physical geometry with some additional
structure, defined on the physical geometry. In the given case the
coordinate system plays the role of this additional structure.

In any case the $\sigma $-Minkowskian geometry, considered as a physical
geometry is a geometry, which is multivariant with respect to spacelike and
null vectors. Thus, the equivalence relation (\ref{b1.7}), (\ref{b1.8}) is
true in the space-time $\sigma $-Minkowskian geometry, and the $\sigma $%
-Minkowskian geometry is multivariant with respect to spacelike and null
vectors.

Working with the space-time geometry of Minkowski, physicists use the
equivalence relation (\ref{b1.18}) and consider the space-time geometry as a
single-variant (transitive) geometry. Such approach appears to be correct,
because one considers only particles with timelike world lines. All vectors,
associated with timelike world lines, are timelike. For these vectors the
equivalence relation (\ref{b1.18}) coincides with the equivalence relation (%
\ref{b1.7}), (\ref{b1.8}).

However, for consideration of free particles, whose world line is a
spacelike helix with the timelike axis, one needs the true equivalence
relation (\ref{b1.7}), (\ref{b1.8}). Strictly, the free particles with such
a helical world line cannot exist in the $\sigma $-Minkowskian space-time.
However, a slight modification of the $\sigma $-Minkowskian space-time
geometry admits existence of free particle with the spacelike helical world
line \cite{R2007c}. Such a particle is associated with the classical Dirac
particle \cite{R2004,R2004a}. At the conventional approach the Dirac
particle is associated with quantum physics, and it has no simple
geometrical interpretation, because the $\gamma $-matrices, which are used
in the description of the Dirac particle, have no direct geometrical
interpretation in terms of points of the space-time.

In general, working with the geometry of Minkowski, one meets the
alternative:

1. The geometry of Minkowski ($\sigma $-Minkowskian geometry) is a physical
geometry, then it is multivariant and it is not axiomatizable.

2. The geometry of Minkowski is single-variant, and axiomatizable, Then it
is not a physical geometry.

Physicists are apt to choose the first possibility, because they need the
geometry for description of the space-time properties, where the main
quantity is the space-time distance (world function). The circumstance,
whether or not the space-time geometry is axiomatizable is a secondary
matter.

On the contrary, mathematicians are apt to choose the second variant. It is
very important for them, that the geometry can be axiomatized, and one can
invent new axioms and prove new theorems. Problems of application of the
geometry in physics and the problem, whether or not the considered geometry
is a physical one, are secondary problems for mathematicians.

It should note, that demanding the geometry axiomatization, mathematicians
have many problems, which could be get round easily. They believe, that any
geometry is axiomatizable, that is it can be deduced from a proper
axiomatics. Such an approach generates the problem of the geometry
consistency and problem of formulation and proof of numerous geometrical
theorems. A use of physical geometries, which are constructed by a
deformation of the proper Euclidean geometry, is free of these problems.
Using the method of the geometry construction, based on the deformation of
the proper Euclidean geometry (the deformation principle), one can construct
physical geometries, which cannot be constructed by means of the
axiomatization method.

\section{Geometry or method of the geometry \newline
description?}

As a rule the mathematicians do not perceive the intransitive equivalence
relation. They say: "In all mathematical papers the equivalence relation is
defined as a transitive operation. If you use intransitive operation, use
another term for this operation, for instance, "generalized equivalence", or
any term other, than the term 'equivalence'. It is necessary, in order to
avoid misunderstanding." The mathematician would be right, if one constructs
another abstract logical construction. However, one discusses properties of
real object of the space-time geometry, and one should use concepts, which
are adequate for description of the space-time. Misunderstanding arises,
because mathematicians consider the linear vector space with the scalar
product, given on it, as a necessary attribute of a physical geometry, and
insist, that all concepts, which are used for description of the linear
vector space, are adequate for description of a physical geometry
(space-time geometry).

In reality, the linear vector space with the scalar product, given on it, is
not an attribute of the physical geometry, in general. Some concepts of the
linear vector space have another meaning, than the meaning of these
concepts, used in dynamics and in the physical geometry. Let us adduce
another examples of arguments against the physical geometry. Vector $\mathbf{%
PQ}$ is a ordered set of to points $\left\{ P,Q\right\} $ in physical
geometry. In conventional approach to geometry the vector is an element of
the linear vector space. As a result mathematicians argue against definition
of a vector as a set of two points.

Let us consider two different definitions of the scalar product of two
vectors $\mathbf{P}_{0}\mathbf{P}_{1}$ and $\mathbf{Q}_{0}\mathbf{Q}_{1}$,
whose coordinate representation has the form (\ref{b1.12}). The conventional
definition, based on the conception of the linear vector space has the form%
\begin{equation}
\left( \mathbf{P}_{0}\mathbf{P}_{1}.\mathbf{Q}_{0}\mathbf{Q}_{1}\right)
=x^{k}y_{k}  \label{b2.0}
\end{equation}

Alternative definition has the form (\ref{b1.5}). The definition (\ref{b2.0}%
) is impossible, if the conception of the linear vector field is not used.
The definition (\ref{b1.5}) is impossible, if the geometry is not a physical
geometry, and the conception of the world function cannot be introduced.
When the geometry is based on the conception of the linear vector field, and
besides, the geometry is physical, both definitions coincide. What of two
definitions is more fundamental? If one considers the geometry as a
space-time geometry (physical geometry), the definition (\ref{b1.5}) is a
true definition of the scalar product. In this case the conventional
definition, based on conception of the linear vector field may not exist,
and one is not to take it into account. If the mathematical geometry is
considered, and the property of the axiomatization is the main property of
the geometry, then the definition (\ref{b2.0}) is the true definition.
However, such a mathematical geometry has no relation to the space-time
geometry of the real world. In any case there is no reason for inventing a
new name for the scalar product (\ref{b1.5}). Like arguments are used in
favour of the definition of a vector, as a ordered set of two points.

The third argument against a physical geometry looks as follows. A segment
of straight has no width (it is one-dimensional). It is one of Euclidean
axioms, and this axiom is essential at the conventional construction of the
linear vector space with the scalar product, given on it.

All arguments against multivariant physical geometry are based on the
supposition, that the linear vector space is a necessary structure of a
physical geometry. In reality the physical geometry is constructed without a
reference to the linear vector space and its properties. There are different
representation of the proper Euclidean geometry \cite{R2007}, and the linear
vector space is only an auxiliary structure, which is used in one of
representations (V-representation) of the proper Euclidean geometry. It
should stress, that the linear vector space with the scalar product, given
on it, \textit{is not a geometry, but only an auxiliary structure}, which is
necessary for introducing the concept of the angle. In other words, the
linear vector space is only one of possible methods, which are used for
description of the proper Euclidean geometry. One may not insist, that
concepts of this description method are adequate in any physical geometry.

Although different representation of the Euclidean geometry are presented in
\cite{R2007}, the principal points of this paper are presented here, because
they are very important for understanding of arguments, adduced by
mathematicians against the intransitive physical geometry.

\section{Three representations of the proper Euclidean geometry}

There are at least three representations of the proper Euclidean geometry,
which differ in the number and in the choice of basic elements (primary
concepts).

The Euclidean representation (E-representation) of the proper Euclidean
geometry contains three basic elements: point, segment and angle. The
segment is a segment of the straight line. It consists of infinite number of
points. The segment is determined uniquely by its end points. The angle is a
figure, formed by two segments provided the end of one segment coincides
with the end of other one. Properties of basic elements are described by a
system of axioms. Any geometrical object may be considered to be some
composition of blocks (\textit{point, segment, angle}). The number of blocks
may be infinite. The segments determine distances. The angles determine the
mutual orientation of segments. Comparison of geometrical objects (figures) $%
\mathcal{O}_{1}$ and $\mathcal{O}_{2}$ is produced by their displacement and
superposition. If two figures coincide at superposition, they are considered
to be equivalent (equal).

The vector representation (V-representation) of the Euclidean geometry
contains two basic elements (\textit{point, vector}). From the viewpoint of
E-representation the vector is a directed segment of straight, determined by
two points. One point is the origin of the vector, another point is the end
of the vector. To construct the angle, there is some auxiliary structure,
known as the linear vector space with the scalar product, given on it. This
structure describes interrelation of two vectors. The mutual orientation of
two vectors $\mathbf{P}_{0}\mathbf{P}_{1}$ and $\mathbf{P}_{0}\mathbf{P}_{2}$
is described, by the angle $\theta $ between them, which is determined by
the relation
\begin{equation}
\left\vert \mathbf{P}_{0}\mathbf{P}_{1}\right\vert \cdot \left\vert \mathbf{P%
}_{0}\mathbf{P}_{2}\right\vert \cos \theta =\left( \mathbf{P}_{0}\mathbf{P}%
_{1}.\mathbf{P}_{0}\mathbf{P}_{2}\right)  \label{b3.1}
\end{equation}%
where $\left( \mathbf{P}_{0}\mathbf{P}_{1}.\mathbf{P}_{0}\mathbf{P}%
_{2}\right) $ is the scalar product of vectors $\mathbf{P}_{0}\mathbf{P}_{1}$
and $\mathbf{P}_{0}\mathbf{P}_{2}$. The quantities $\left\vert \mathbf{P}_{0}%
\mathbf{P}_{1}\right\vert \ $and $\left\vert \mathbf{P}_{0}\mathbf{P}%
_{2}\right\vert $ are their lengths%
\begin{equation}
\left\vert \mathbf{P}_{0}\mathbf{P}_{1}\right\vert ^{2}=\left( \mathbf{P}_{0}%
\mathbf{P}_{1}.\mathbf{P}_{0}\mathbf{P}_{1}\right) ,\qquad \left\vert
\mathbf{P}_{0}\mathbf{P}_{2}\right\vert ^{2}=\left( \mathbf{P}_{0}\mathbf{P}%
_{2}.\mathbf{P}_{0}\mathbf{P}_{2}\right)  \label{b3.2}
\end{equation}

Thus, the transition from the E-representation with three basic elements to
the V-representation with two basic elements is possible, provided the
properties of a vector are determined by the fact, that the vector is an
element of the linear vector space with the scalar product, given on it,
i.e. there is an auxiliary structure, describing interrelation of two
vectors.

It is possible a further reduction of the number of basic elements in the
description of the proper Euclidean geometry. The $\sigma $-representation
(in terms world function) of the proper Euclidean geometry contains only one
basic element (point). Besides, there is a structure, which admits one to
obtain two other basic elements (vector and angle). This structure is the
world function, which describes interrelation of two points.

The transition from the E-representation to the V-representation reduces the
number of basic elements. Simultaneously this transition generates a new
structure (linear vector space). The linear vector space with the scalar
product, given on it, determines those properties of vector, which admits
one to restore the removed angle and its properties.

The transition from the V-representation to the $\sigma $-representation
reduces also the number of basic elements. Only one basic element (\textit{%
point}) remains. Simultaneously this transition replaces the linear vector
space by a new two-point structure (\textit{world function}), which
describes interrelation of two points (instead of two vectors). The $\sigma $%
-representation was not known for a long time. It appeared only in the end
of the twentieth century \cite{R2001}, and its properties are not yet
investigated completely. In the twentieth century mathematicians used the
V-representation as the most developed description of the geometry.

The proper Euclidean geometry $\mathcal{G}_{\mathrm{E}}$ is axiomatizable.
Axiomatics is to depend on the representation, which is used for description
of $\mathcal{G}_{\mathrm{E}}$. In E-representation one has axiomatics $%
\mathcal{A}_{\mathrm{E}}$. In V-representation one has axiomatics $\mathcal{A%
}_{\mathrm{V}}$. In $\sigma $-representation one has axiomatics $\mathcal{A}%
_{\mathrm{\sigma }}$. Axiomatics $\mathcal{A}_{\mathrm{E}},\mathcal{A}_{%
\mathrm{V}}$ and $\mathcal{A}_{\mathrm{\sigma }}$ contain different axioms,
but these axiomatics are equivalent in the sense, that any two axiomatics
can be deduced from the third axiomatics. However, it is true only for the
proper Euclidean geometry. If one changes slightly one of axiomatics, for
instance, $\mathcal{A}_{\mathrm{V}}$, hoping to obtain another geometry (for
instance, pseudo-Euclidean), the axiomatics $\mathcal{A}_{\mathrm{E}},%
\mathcal{A}_{\mathrm{V}},$ $\mathcal{A}_{\mathrm{\sigma }}$ cease to be
equivalent, in general, because the system of axioms (axiomatics) is a
finite set of some propositions of the axiomatizable geometry. Changing in
the first section the axiomatics $\mathcal{A}_{\mathrm{V}}$, to obtain the
pseudo-Euclidean geometry, one discovers, that a simple refuse from
non-negativity of the vector length generates the alternative: either
non-physical axiomatizable geometry, or multivariant physical geometry with
intransitive equality relation. At any rate, after such a very simple
modification one cannot consider axiomatics $\mathcal{A}_{\mathrm{V}}$ and $%
\mathcal{A}_{\mathrm{\sigma }}$ as equivalent.

Basic elements (point, segment, angle) are not independent. Constructing a
generalized geometry in V-representation, one should modify properties of
all basic elements coordinated. In $\sigma $-representation there is only
one basic element (point). The problem of coordinated modification of
different basic elements of the geometry is absent, because there is only
one basic element. Interrelation of two points in the proper Euclidean
geometry may be changed by an arbitrary change of the world function. As far
as all propositions of the proper Euclidean geometry are expressed via the
world function, any change of the world function generates concerted change
of properties of other geometrical objects (including such important objects
as the straight and the equivalence relation).

For instance, in E-representation the straight has no width, and it is an
axiom, which is used in E-representation, as well as in V-representation. In
the $\sigma $-representation the segment $\mathcal{T}_{\left[ P_{0}P_{1}%
\right] }$ of straight between points $P_{0}$ and $P_{1}$ is one of
propositions of the physical geometry, which does not pretend to be an
axiom. This proposition is valid for all physical geometry, including the
proper Euclidean geometry. In the $\sigma $-representation the segment $%
\mathcal{T}_{\left[ P_{0}P_{1}\right] }$ is defined as a set of point $R$%
\begin{equation}
\mathcal{T}_{\left[ P_{0}P_{1}\right] }=\left\{ R|\sqrt{2\sigma \left(
P_{0},R\right) }+\sqrt{2\sigma \left( R,P_{1}\right) }=\sqrt{2\sigma \left(
P_{0},P_{1}\right) }\right\}  \label{b3.3}
\end{equation}%
The definition (\ref{b3.3}) of the segment $\mathcal{T}_{\left[ P_{0}P_{1}%
\right] }$ contains only points $P_{0},P_{1}$, determining the segment, and
the world function between the points $P_{0},P_{1}$ and the running point $R$
of the set $\mathcal{T}_{\left[ P_{0}P_{1}\right] }$.

In general, the set (\ref{b3.3}) is a surface. Dimension of this surface
depends on the dimension of the set $\Omega $ and on the form of the world
function. If $\sigma $ is the world function of the proper Euclidean
geometry ($\sigma =\sigma _{\mathrm{E}}$), the straight (\ref{b3.3}) is
one-dimensional line independently of the dimension of $\Omega $. This
property is the special property of the proper Euclidean geometry, which is
determined by the properties of the Euclidean world function $\sigma _{%
\mathrm{E}}$. If the world function $\sigma \neq \sigma _{\mathrm{E}}$, the
set (\ref{b3.3}) is not one-dimensional, in general. In the case, when $%
\Omega $ is the $n$-dimensional manifold, the dimension of the set (\ref%
{b3.3}) is equal to $n-1$, in general. The axiom of axiomatics $\mathcal{A}_{%
\mathrm{E}}$ and $\mathcal{A}_{\mathrm{V}}$, which declares, that the
straight has no width, is incompatible with (\ref{b3.3}), in general. The
axiom on one-dimensionality of the set $\mathcal{T}_{\left[ P_{0}P_{1}\right]
}$ is not valid, in general, in the deformed physical geometry, described by
the world function $\sigma \neq \sigma _{\mathrm{E}}$.

From viewpoint of the physical geometry one-dimensionality of the segment (%
\ref{b3.3}) is a corollary of the definition (\ref{b3.3}) and of the special
properties of the Euclidean world function $\sigma _{\mathrm{E}}$. If the
world function $\sigma _{\mathrm{E}}$ is replaced by another world function $%
\sigma $, the special properties of $\sigma _{\mathrm{E}}$ are violated, in
general, and the segment $\mathcal{T}_{\left[ P_{0}P_{1}\right] }$ ceases to
be one-dimensional.

\section{Peculiarities and capacities of multivariant \newline
physical geometries}

The multivariant physical geometries, i.e. physical geometries with
intransitive equivalence relation have many such capacities, which cannot be
possessed by axiomatizable geometries. The multivariant physical geometry
cannot be an axiomatizable geometry, because in any axiomatizable geometry
the equivalence relation is transitive, and the axiomatizable geometry is
single-variant. The class of possible physical (multivariant) geometries is
more powerful, than the class of axiomatizable geometries. For instance, the
class of homogeneous isotropic space-time geometries is labelled by a
function of one argument. Any homogeneous isotropic space-time geometry is
described by the world function of the form%
\begin{equation}
\sigma =F\left( \sigma _{\mathrm{M}}\right)  \label{b4.1}
\end{equation}%
where $\sigma _{\mathrm{M}}$ is the world function of the geometry of
Minkowski, and $F$ is an arbitrary real function, which has the property $%
F\left( 0\right) =0$. Among these geometries there is only one Riemannian
geometry -- the geometry of Minkowski. Other physical geometries are not
Riemannian. As we have seen in the first section, the geometry of Minkowski
is either physical, then it is multivariant and cannot be axiomatizable, or
it is axiomatizable, but then it is not a physical geometry.

The multivariant space-time geometry may be discrete, and this discrete
geometry may be given on the continual manifold of Minkowski. World function
$\sigma _{\mathrm{d}}$ of the discrete space-time geometry $\mathcal{G}_{%
\mathrm{d}}$ is expressed via the Minkowskian world function $\sigma _{%
\mathrm{M}}$. For instance, the space-time geometry with the world function
\begin{equation}
\sigma _{\mathrm{d}}=\sigma _{\mathrm{M}}+d\mathrm{sgn}\left( \sigma _{%
\mathrm{M}}\right) ,\qquad d\equiv \lambda _{0}^{2}=\frac{\hbar }{2bc}=\text{%
const}  \label{b4.2}
\end{equation}%
\begin{equation}
\text{sgn}\left( x\right) =\left\{
\begin{array}{lll}
\frac{x}{\left\vert x\right\vert } & \text{if} & x\neq 0 \\
0 & \text{if} & x=0%
\end{array}%
\right.  \label{b4.3}
\end{equation}%
is discrete, because in this geometry there are no points, separated by the
distance, which is less, than $\sqrt{2}\lambda _{0}$. This world function $%
\sigma _{\mathrm{d}}$ is Lorentz-invariant, as far as $\sigma _{\mathrm{d}}$
is a function of $\sigma _{\mathrm{M}}$, and $\sigma _{\mathrm{M}}$ is
Lorentz-invariant. Usually it is supposed that a discrete geometry is a
geometry, given on the point set of the type of a grid, what is incompatible
with the Lorentz-invariance. Discreteness of the space-time geometry $%
\mathcal{G}_{\mathrm{d}}$ seems to be very surprising, because the
space-time is given on the continual manifold of Minkowski. Conventionally
the discrete space is associated with a grid. The discrete space-time on the
continuous manifold seems to be impossible. This example shows, that a
physical geometry and a continuous manifold, where the geometry is given,
are quite different matters. Manifold and its dimension are only attributes
of the vector representation of the Euclidean geometry (i.e. attributes of
the description method), whereas the discreteness of a geometry is an
attribute of the geometry in itself. Besides, motion of free particles in
the space-time geometry (\ref{b4.2}) appears to be multivariant
(stochastic), and a statistical description of this multivariant motion is
equivalent to quantum description in terms of the Schr\"{o}dinger equation
\cite{R91}. The space-time geometry (\ref{b4.2}) contains the elementary
length $\lambda _{0}$, which is expressed via the quantum constant $\hbar $,
via the speed of the light $c$ and via the universal constant $b$. The
universal constant $b$ connects the geometrical length $\mu $ of the
momentum vector $\mathbf{P}_{k}\mathbf{P}_{k+1}$ with the mass $m$,
associated with momentum vector $\mathbf{P}_{k}\mathbf{P}_{k+1}$, by means
of the relation
\begin{equation}
m=b\mu =b\left\vert \mathbf{P}_{k}\mathbf{P}_{k+1}\right\vert  \label{b4.4}
\end{equation}%
Apparently, the quantum effects in the motion of free particles and
discreteness of the space-time geometry are connected matters.

Discreteness of the space-time geometry may be incomplete. For instance, let
us consider the space-time geometry $\mathcal{G}_{\mathrm{g}}$, described by
the world function%
\begin{equation}
\sigma _{\mathrm{g}}=\sigma _{\mathrm{M}}+\lambda _{0}^{2}\left\{
\begin{array}{lll}
\mathrm{sgn}\left( \sigma _{\mathrm{M}}\right) & \text{if} & \left\vert
\sigma _{\mathrm{M}}\right\vert >\sigma _{0} \\
\frac{\sigma _{\mathrm{M}}}{\sigma _{0}} & \text{if} & \left\vert \sigma _{%
\mathrm{M}}\right\vert \leq \sigma _{0}%
\end{array}%
\right. ,\qquad \lambda _{0}^{2},\sigma _{0}=\text{const}\geq 0  \label{b4.5}
\end{equation}%
If the constant $\sigma _{0}\rightarrow 0$, the space-time geometry (\ref%
{b4.5}) tends to the discrete geometry (\ref{b4.2}). As it is follows from (%
\ref{b4.5}), the relative density $\rho \left( \sigma _{\mathrm{g}}\right) $
of points in the geometry $\mathcal{G}_{\mathrm{g}}$ with respect to the
geometry of Minkowski $\mathcal{G}_{\mathrm{M}}$ is described by the relation%
\begin{equation}
\rho \left( \sigma _{\mathrm{g}}\right) =\frac{d\sigma _{\mathrm{M}}\left(
\sigma _{\mathrm{g}}\right) }{d\sigma _{\mathrm{g}}}=\left\{
\begin{array}{lll}
1 & \text{if} & \left\vert \sigma _{\mathrm{g}}\right\vert >\sigma
_{0}+\lambda _{0}^{2} \\
\frac{\sigma _{0}}{\sigma _{0}+\lambda _{0}^{2}} & \text{if} & \left\vert
\sigma _{\mathrm{g}}\right\vert \leq \sigma _{0}+\lambda _{0}^{2}%
\end{array}%
\right.  \label{b4.6}
\end{equation}%
One can see from the relation (\ref{b4.6}), that the relative density $\rho
\left( \sigma _{\mathrm{g}}\right) $ of particles for $\sigma _{\mathrm{g}%
}\in \left( -\sigma _{0}-\lambda _{0}^{2},\sigma _{0}+\lambda
_{0}^{2}\right) $ is less, than unity. The relative density $\rho \left(
\sigma _{\mathrm{g}}\right) \rightarrow 0$, when $\sigma _{0}\rightarrow 0$,
and the geometry $\mathcal{G}_{\mathrm{g}}$ tends to the discrete geometry (%
\ref{b4.2}). The space-time geometry $\mathcal{G}_{\mathrm{g}}$, where the
relative density $\rho \left( \sigma _{\mathrm{g}}\right) $ satisfies the
condition $\left\vert \rho \left( \sigma _{\mathrm{g}}\right) \right\vert <1$
for $\left\vert \sigma _{\mathrm{g}}\right\vert <\lambda _{0}^{2}$, should
be qualified as a partly discrete, or grainy space-time geometry. The grainy
space-time geometry cannot be realized in terms of the geometry, based on
the linear vector space.

Nevertheless, in the framework of the multivariant physical geometries the
grainy geometry is a usual space-time geometry, where all concepts and
objects of the proper Euclidean geometry may be introduced. It is quite
another matter, that some of these objects look rather exotic and unusual.

The pseudo-Riemannian space-time geometry is multivariant with respect to
spacelike vectors, as we have seen in the example of the $\sigma $%
-Minkowskian geometry. However, the proper Riemannian physical geometry is
also multivariant,in general, because of its curvature. In general, the
multivariance of a physical geometry manifests itself in
non-one-dimensionality of straights. The straight (geodesic) $\mathcal{T}%
_{P_{0}P_{1}}$, passing through points $P_{0}$ and $P_{1}$ is
one-dimensional in the proper Riemannian geometry. The one-dimensionality of
the geodesic is conditioned by the special property of the Riemannian world
function%
\begin{equation}
\sigma _{\mathrm{R}}\left( P_{0},P_{1}\right) =\frac{1}{2}\left(
\dint\limits_{\mathcal{T}_{\left[ P_{0}P_{1}\right] }}\sqrt{g_{ik}\left(
x\right) dx^{i}dx^{k}}\right) ^{2}  \label{b4.7}
\end{equation}%
where integration is produced along the geodesic segment $\mathcal{T}_{\left[
P_{0}P_{1}\right] }$, connecting points $P_{0}$ and $P_{1}$.

The world function, defined by the relation (\ref{b4.7}) satisfies the
triangle axiom%
\begin{equation}
\sqrt{2\sigma _{\mathrm{R}}\left( P_{0},R\right) }+\sqrt{2\sigma _{\mathrm{R}%
}\left( R,P_{1}\right) }\geq \sqrt{2\sigma _{\mathrm{R}}\left(
P_{0},P_{1}\right) },\qquad \forall P_{0},P_{1},R\in \Omega  \label{b4.8}
\end{equation}%
and this property generates the one-dimensionality of the straight $\mathcal{%
T}_{P_{0};P_{0}P_{1}}$, passing through the point $P_{0}$ collinear to
vector $\mathbf{P}_{0}\mathbf{P}_{1}$. In the general case the straight
passing through the point $Q_{0}$, collinear to vector $\mathbf{P}_{0}%
\mathbf{P}_{1}$ is defined by the relation
\begin{equation}
\mathcal{T}_{Q_{0};P_{0}P_{1}}=\left\{ R|\mathbf{Q}_{0}\mathbf{R}\parallel
\mathbf{P}_{0}\mathbf{P}_{1}\right\}  \label{b4.9}
\end{equation}%
where the collinearity condition $\mathbf{Q}_{0}\mathbf{R}\parallel \mathbf{P%
}_{0}\mathbf{P}_{1}$ is described by the relation%
\begin{eqnarray}
\mathbf{Q}_{0}\mathbf{R} &\parallel &\mathbf{P}_{0}\mathbf{P}_{1}\mathbf{%
:\qquad }\left\vert
\begin{array}{ll}
\left( \mathbf{Q}_{0}\mathbf{R}.\mathbf{Q}_{0}\mathbf{R}\right) & \left(
\mathbf{Q}_{0}\mathbf{R}.\mathbf{P}_{0}\mathbf{P}_{1}\right) \\
\left( \mathbf{P}_{0}\mathbf{P}_{1}.\mathbf{Q}_{0}\mathbf{R}\right) & \left(
\mathbf{P}_{0}\mathbf{P}_{1}.\mathbf{P}_{0}\mathbf{P}_{1}\right)%
\end{array}%
\right\vert  \notag \\
&=&\left\vert \mathbf{P}_{0}\mathbf{P}_{1}\right\vert ^{2}\left\vert \mathbf{%
Q}_{0}\mathbf{R}\right\vert ^{2}-\left( \mathbf{Q}_{0}\mathbf{R}.\mathbf{P}%
_{0}\mathbf{P}_{1}\right) ^{2}=0  \label{b4.10}
\end{eqnarray}%
In the Riemannian geometry only the straight $\mathcal{T}_{P_{0};P_{0}P_{1}}$
is one-dimensional, whereas the straight (\ref{b4.9}) with $P_{0}\neq Q_{0}$
is not one-dimensional, in general. This circumstance is a manifestation of
the multivariance of the Riemannian geometry, considered as a physical
geometry. However, the Riemannian geometry is not considered usually as a
physical geometry with the world function, defined by the relation (\ref%
{b4.7}). To avoid a confusion, we shall use the term "$\sigma $-Riemannian
geometry" for the physical geometry with the world function (\ref{b4.7}),
remaining the term "Riemannian geometry" for the logical construction in its
conventional form, when one introduces infinitesimal interval on a manifold.
Many propositions of the $\sigma $-Riemannian geometry and those of the
Riemannian geometry, concerning dynamics in the space-time, coincide. They
differ in the field, concerning fernparallelism (parallelism of remote
vectors). In the $\sigma $-Riemannian geometry there is the absolute
parallelism, and it is multivariant. In the Riemannian geometry the
parallelism of remote vectors is forbidden (it is simply not considered).
Then the straight (\ref{b4.9}) with $P_{0}\neq Q_{0}$ cannot by constructed.
Instead of the absolute parallelism, in the Riemannian geometry the parallel
transport is introduced. Result of the parallel transport depends on the
path of the transport, i.e. essentially the parallelism of remote vectors is
multivariant. However, one cannot accept multivariance of the parallelism,
because the construction of the Riemannian geometry is based on the
transitive equivalence relation. Absolute parallelism, which would be
multivariant with necessity, is incompatible with the transitivity of the
equivalence relation. However, the prohibition of the absolute parallelism
(non-consideration of parallelism of remote vectors) does not eliminate
inconsistency between the multivariance of parallelism and the transitive
equivalence relation, because the inconsistency does not disappear by the
reason of the fact, that one does not consider corollaries of this
inconsistency.

The $\sigma $-Riemannian geometry has no problems with inconsistency,
because it is a physical geometry, which is multivariant and which is not
axiomatizable. The physical geometry, which is not deduced from axioms,
cannot be inconsistent in principle, because inconsistency means, that two
different methods of deduction of some proposition lead to different results.

In any mathematical model the equivalence relation is transitive. This
property of the mathematical model provides definiteness (single-variance)
for all conclusions, made on the basis of such a mathematical model. If the
equivalence relation is intransitive, the conclusions, made on the basis of
such a model, cease to be definite. They becomes multivariant. The logical
construction with the intransitive equivalence relation and, hence, with
multivariant conclusions is not considered to be a mathematical model,
because it is non-causal and useless. One cannot make a definite prediction
on the basis of such a model. Besides, such a model cannot be axiomatizable,
because the axiomatization supposes a single-variance of conclusions.

Fortunately, a multivariant model can be reduced to a single-variant model,
provided one unites the set of many conclusions, which follows from one
statement into one conclusion. In other words, one considers the set of
different objects as an united object (statistical ensemble). One may work
with the statistical ensemble, considering it as a single object. Then the
multivariant model may cease to be multivariant. It turns into a
single-variant (transitive) model provided, that its objects be statistical
ensembles of original objects. Such a procedure is known as the statistical
description, which deals with statistically averaged objects. Prediction of
the model about statistically averaged objects (statistical ensembles),
which are now objects of the model, may appear to be single-variant, if the
statistical description is produced properly. In other words, a statistical
description, produced properly, transforms a multivariant model into a
single-variant mathematical model. Procedure of the statistical description
is well known. It is used in different branches of theoretical physics.

In the physical geometry along with the multivariance of the equivalence of
two vectors, there may exist a zero-variance, when at the point $P_{0}$
there is no vector, which is equivalent to a given vector $\mathbf{Q}_{0}%
\mathbf{Q}_{1}$. If multivariance may be taken in to account in dynamics be
means of a proper statistical description, the zero-variance cannot be taken
into account in this way. The zero-variance generated some discrimination
mechanism, which prohibits some values of dynamic parameters (masses,
charges, angular momenta) \cite{R2007c,R2008}. Statistical description can
be taken into account in framework of quantum principles and quantum
mechanics, whereas the discrimination mechanism cannot be taken into account
by the quantum principles.

To understand, why the discrimination mechanism is something external with
respect to quantum dynamics, one should keep in mind, that the dynamics in
the physical space-time geometry distinguishes from the conventional
dynamics in the space-time, constructed in the framework of the
V-representation. The conventional dynamics is founded on the conception of
the linear vector space. It supposes unlimited divisibility of the
space-time, which is associated with a description in terms of differential
equations.

In the space-time, described by a physical geometry, which may be grainy and
have a limited divisibility, the dynamics is formulated originally in terms
of finite (not differential) dynamic equations. For instance, if the
space-time is discrete, it would be strange to demand, that dynamic equation
be formulated in terms of differential equations, because the infinitesimal
distances are absent in the discrete space-time.

\section{Dynamics in the physical space-time geometry}

In the space-time geometry of Minkowski a motion of a pointlike particle is
described by its world line $\mathcal{L}$, which may be considered as a
chain of infinitesimal connected vectors $\mathbf{P}_{k}\mathbf{P}_{k+1}$, $%
k=..0,1,...$. All vectors $\mathbf{P}_{k}\mathbf{P}_{k+1}$ are infinitesimal
and timelike. If the motion of the particle is free, the adjacent vectors $%
\mathbf{P}_{k}\mathbf{P}_{k+1}$ are equivalent. As far as infinitesimal
vectors $\mathbf{P}_{k}\mathbf{P}_{k+1}$ are timelike, their equivalence is
single-variant. Besides, some non-geometrical quantity $m$, known as the
particle mass, is ascribed to the pointlike particle. The vector $\mathbf{P}%
_{k}\mathbf{P}_{k+1}$ describes the state of the pointlike particle at the
point $P_{k}\in \mathcal{L}$. As far as the vector $\mathbf{P}_{k}\mathbf{P}%
_{k+1}$ is infinitesimal, one uses a finite vector $\mathbf{p}$, which is in
parallel with the infinitesimal vector $\mathbf{P}_{k}\mathbf{P}_{k+1}$ and
has the length equal to the mass $m$. Vector $\mathbf{p}$ is known as the
momentum 4-vector.

In the physical space-time geometry a geometrical object is a subset of
points of the space-time $\Omega $. Geometrical object is considered as a
composition of elementary objects (EGO). Any elementary geometrical object $%
\mathcal{E}$ is defined by its skeleton $\mathcal{P}^{n}\equiv \left\{
P_{0},P_{1},...P_{n}\right\} \subset \Omega $ and its envelope function $f_{%
\mathcal{P}^{n}}$
\begin{equation}
f_{\mathcal{P}^{n}}:\qquad \Omega \rightarrow \mathbb{R}  \label{b5.1}
\end{equation}%
The skeleton $\mathcal{P}^{n}$ is the set of $n+1$ points of the point set $%
\Omega $. The envelope function is a function of the running point $R\in
\Omega $ and of parameters $\mathcal{P}^{n}\subset \Omega $. The envelope
function $f_{\mathcal{P}^{n}}$ is supposed to be an algebraic function of $s$
arguments $w=\left\{ w_{1},w_{2},...w_{s}\right\} $, $s=(n+2)(n+1)/2$. Each
of arguments $w_{k}=\sigma \left( Q_{k},L_{k}\right) $ is a world function
of two arguments $Q_{k},L_{k}\in \left\{ R,\mathcal{P}^{n}\right\} $, either
belonging to skeleton $\mathcal{P}^{n}$, or coinciding with the running
point $R$. Thus, any elementary geometric object $\mathcal{E}$ is determined
by its skeleton and its envelope function as the set of zeros of the
envelope function
\begin{equation}
\mathcal{E}=\left\{ R|f_{\mathcal{P}^{n}}\left( R\right) =0\right\} ,\qquad
,R\in \Omega ,\quad \mathcal{E\subset }\Omega ,  \label{b5.2}
\end{equation}

For instance, the circular cylinder $\mathcal{C}(P_{0},P_{1},Q)$ with the
points $P_{0},P_{1}$ on the cylinder axis and the point $Q$ on its surface
is determined by the relation
\begin{eqnarray}
\mathcal{C}(P_{0},P_{1},Q) &=&\left\{ R|f_{P_{0}P_{1}Q}\left( R\right)
=0\right\} ,\qquad \mathcal{C}(P_{0},P_{1},Q)\subset \Omega  \label{g3.1} \\
f_{P_{0}P_{1}Q}\left( R\right) &=&F_{2}\left( P_{0},P_{1},Q\right)
-F_{2}\left( P_{0},P_{1},R\right)  \notag
\end{eqnarray}%
\begin{equation}
F_{2}\left( P_{0},P_{1},Q\right) =\left\vert
\begin{array}{cc}
\left( \mathbf{P}_{0}\mathbf{P}_{1}.\mathbf{P}_{0}\mathbf{P}_{1}\right) &
\left( \mathbf{P}_{0}\mathbf{P}_{1}.\mathbf{P}_{0}\mathbf{Q}\right) \\
\left( \mathbf{P}_{0}\mathbf{Q}.\mathbf{P}_{0}\mathbf{P}_{1}\right) & \left(
\mathbf{P}_{0}\mathbf{Q}.\mathbf{P}_{0}\mathbf{Q}\right)%
\end{array}%
\right\vert  \label{g3.2}
\end{equation}%
Here $\sqrt{F_{2}\left( P_{0},P_{1},Q\right) }$ is the area of the
parallelogram, constructed on the vectors $\mathbf{P}_{0}\mathbf{P}_{1}$, $%
\mathbf{P}_{0}\mathbf{Q}$,\textbf{\ }and $\frac{1}{2}\sqrt{F_{2}\left(
P_{0},P_{1},Q\right) }\ $ is the area of triangle with vertices at the
points $P_{0},P_{1},Q$. The equality $F_{2}\left( P_{0},P_{1},Q\right)
=F_{2}\left( P_{0},P_{1},R\right) $ means that the distance between the
point $Q$ and the axis, determined by the vector $\mathbf{P}_{0}\mathbf{P}%
_{1}$ is equal to the distance between $R$ and the axis.

The elementary geometrical object $\mathcal{E}$ is determined in all
physical geometries at once. In particular, it is determined in the proper
Euclidean geometry, where one can obtain its meaning. One interprets the
elementary geometrical object $\mathcal{E}$, using our knowledge of the
proper Euclidean geometry. Thus, the proper Euclidean geometry is used as a
sample geometry for interpretation of any physical geometry. In particular,
the cylinder (\ref{g3.1}) is determined uniquely in any physical geometry
with any world function $\sigma .$

In the Euclidean geometry the points $P_{0}$ and $P_{1}$ determine the
cylinder axis. In the proper Euclidean geometry the shape of a cylinder
depends on its axis and radius, but not on the disposition of points $%
P_{0},P_{1}$ on the cylinder axis. As a result in the Euclidean geometry the
cylinders $\mathcal{C}(P_{0},P_{1},Q)$ and $\mathcal{C}(P_{0},P_{2},Q)$
coincide, provided vectors $\mathbf{P}_{0}\mathbf{P}_{1}$ and $\mathbf{P}_{0}%
\mathbf{P}_{2}$ are collinear. In the general case of physical geometry the
cylinders $\mathcal{C}(P_{0},P_{1},Q)$ and $\mathcal{C}(P_{0},P_{2},Q)$ do
not coincide, in general, even if vectors $\mathbf{P}_{0}\mathbf{P}_{1}$ and
$\mathbf{P}_{0}\mathbf{P}_{2}$ are collinear. Thus, in general, a
deformation of the Euclidean geometry splits Euclidean geometrical objects.
It is a manifestation of the multivariance of a physical geometry.

Two elementary geometrical objects $\mathcal{E}_{\mathcal{P}^{n}}^{\left(
1\right) }$ and $\mathcal{E}_{\mathcal{Q}^{n}}^{\left( 2\right) }$ are
equivalent, if their skeletons $\mathcal{P}^{n}=\left\{
P_{0},P_{1},...P_{n}\right\} $ and $\mathcal{Q}^{n}=\left\{
Q_{0},Q_{1},...Q_{n}\right\} $ are equivalent%
\begin{equation}
\mathcal{P}^{n}\mathrm{eqv}\mathcal{Q}^{n}:\qquad \mathbf{P}_{i}\mathbf{P}%
_{k}\mathrm{eqv}\mathbf{Q}_{i}\mathbf{Q}_{k},\qquad i,k=0,1,...n
\label{b5.3}
\end{equation}%
and their envelope functions $f_{\mathcal{P}^{n}}^{\left( 1\right) }$ and $%
f_{\mathcal{Q}^{n}}^{\left( 2\right) }$ are equivalent
\begin{equation}
f_{\mathcal{P}^{n}}^{\left( 1\right) }\left( R\right) =\Phi \left( f_{%
\mathcal{Q}^{n}}^{\left( 2\right) }\left( R\right) \right) ,\qquad \Phi
\left( 0\right) =0,\qquad \forall R\in \Omega  \label{b5.4}
\end{equation}

Only equivalence of skeletons is important in the formulation of geometrical
objects dynamics, and we shall not distinguish between the geometrical
object and its skeleton. The evolution of a free geometrical object is
described by a chain $\mathcal{C}$ of connected geometrical objects
(skeletons) ...$\mathcal{P}_{\left( 1\right) }^{n},\mathcal{P}_{\left(
2\right) }^{n},...\mathcal{P}_{\left( s\right) }^{n},....$The links
(geometical objects) of the world chain $\mathcal{C}$ are connected in the
sense, that%
\begin{equation}
P_{1}^{\left( s\right) }=P_{0}^{\left( s+1\right) },\qquad s=...0,1,...
\label{b5.5}
\end{equation}%
The vector $\mathbf{P}_{0}^{\left( s\right) }\mathbf{P}_{1}^{\left( s\right)
}=\mathbf{P}_{0}^{\left( s\right) }\mathbf{P}_{0}^{\left( s+1\right) }$,
connecting adjacent links of the chain, is called the leading vector. As far
as the physical space-time geometry is multivariant, in general, only
adjacent links of the world chain are equivalent. The geometrical object,
whose skeleton contains more, than two points will be interpreted as a
composite particle. A pointlike particle is described by two points $\mathbf{%
P}_{s}\mathbf{P}_{s+1}=\left\{ P_{s},P_{s+1}\right\} $. The vector $\mathbf{P%
}_{s}\mathbf{P}_{s+1}$ is interpreted as the geometrical momentum of the
particle, and $\mu =\left\vert \mathbf{P}_{s}\mathbf{P}_{s+1}\right\vert $
is the geometrical mass of the pointlike particle. The geometrical mass $\mu
$ is connected with the usual mass $m$ of the pointlike particle by the
relation (\ref{b4.4}).

Dynamics of free composite particles is formulated in terms of only
geometrical concepts (points and world functions between them). Any
reference to means of description (coordinates, dimension) is absent. As a
result any consideration of invariance with respect to coordinate
transformations and with respect any mappings has no relation to dynamics.
Dynamics of a free composite particle, described by the skeleton $\mathcal{P}%
^{n}$, is a translation without rotation of all skeleton points at the
distance, described by the leading vector $\mathbf{P}_{0}\mathbf{P}_{1}$.
The rotation is absent only from the viewpoint of the geometry $\mathcal{G}$%
, described by the world function $\sigma $. The same world chain, described
in terms of the geometry of Minkowski $\mathcal{G}_{\mathrm{M}}$, contains
translation accompanied by the rotation.

Dynamics of free composite particle in a physical space-time geometry can be
reduced to dynamics of a composite particle, which moves in space-time
geometry of Minkowski under some metrical forces. Let the real world
function has the form
\begin{equation}
\sigma \left( P_{0},P_{1}\right) =\sigma _{\mathrm{M}}\left(
P_{0},P_{1}\right) +d\left( P_{0},P_{1}\right) ,\qquad \forall
P_{0},P_{1}\in \Omega  \label{b5.6}
\end{equation}%
The condition of equivalence (\ref{b5.3}) can be written in the form%
\begin{eqnarray}
\left( \mathbf{P}_{k}^{\left( s\right) }\mathbf{P}_{l}^{\left( s\right) }.%
\mathbf{P}_{k}^{\left( s+1\right) }\mathbf{P}_{l}^{\left( s+1\right)
}\right) &=&\left\vert \mathbf{P}_{k}^{\left( s\right) }\mathbf{P}%
_{l}^{\left( s\right) }\right\vert ^{2},\qquad k,l=0,1,...n,\qquad
s=...0,1,...  \label{b5.6a} \\
\left\vert \mathbf{P}_{k}^{\left( s\right) }\mathbf{P}_{l}^{\left( s\right)
}\right\vert ^{2} &=&\left\vert \mathbf{P}_{k}^{\left( s+1\right) }\mathbf{P}%
_{l}^{\left( s+1\right) }\right\vert ^{2},\quad k,l=0,1,...n,\quad
s=...0,1,...  \label{b5.6b}
\end{eqnarray}%
Using (\ref{b5.6}), one can rewrite the relations (\ref{b5.6a}), (\ref{b5.6b}%
) in the form%
\begin{eqnarray}
&&\left( \mathbf{P}_{k}^{\left( s\right) }\mathbf{P}_{l}^{\left( s\right) }.%
\mathbf{P}_{k}^{\left( s+1\right) }\mathbf{P}_{l}^{\left( s+1\right)
}\right) _{\mathrm{M}}+w\left( P_{k}^{\left( s\right) },P_{l}^{\left(
s\right) },P_{k}^{\left( s+1\right) },P_{l}^{\left( s+1\right) }\right)
\notag \\
&=&2\sigma _{\mathrm{M}}\left( P_{k}^{\left( s\right) },P_{l}^{\left(
s\right) }\right) +2d\left( P_{k}^{\left( s\right) },P_{l}^{\left( s\right)
}\right)  \label{b5.7}
\end{eqnarray}%
\begin{equation}
2\sigma _{\mathrm{M}}\left( P_{k}^{\left( s\right) },P_{l}^{\left( s\right)
}\right) +2d\left( P_{k}^{\left( s\right) },P_{l}^{\left( s\right) }\right)
=2\sigma _{\mathrm{M}}\left( P_{k}^{\left( s+1\right) },P_{l}^{\left(
s+1\right) }\right) +2d\left( P_{k}^{\left( s+1\right) },P_{l}^{\left(
s+1\right) }\right)  \label{b5.8}
\end{equation}%
where%
\begin{eqnarray}
w\left( P_{k}^{\left( s\right) },P_{l}^{\left( s\right) },P_{k}^{\left(
s+1\right) },P_{l}^{\left( s+1\right) }\right) &=&d\left( P_{k}^{\left(
s\right) },P_{l}^{\left( s+1\right) }\right) +d\left( P_{l}^{\left( s\right)
},P_{k}^{\left( s+1\right) }\right)  \notag \\
&&-d\left( P_{k}^{\left( s\right) },P_{k}^{\left( s+1\right) }\right)
-d\left( P_{l}^{\left( s\right) },P_{l}^{\left( s+1\right) }\right)
\label{b5.9}
\end{eqnarray}%
Index "M" at the scalar product means that the scalar product is calculated
in the geometry of Minkowski. The scalar product $\left( \mathbf{P}%
_{k}^{\left( s\right) }\mathbf{P}_{l}^{\left( s\right) }.\mathbf{P}%
_{k}^{\left( s+1\right) }\mathbf{P}_{l}^{\left( s+1\right) }\right) _{%
\mathrm{M}}$ can be represented in the form%
\begin{equation}
\left( \mathbf{P}_{k}^{\left( s\right) }\mathbf{P}_{l}^{\left( s\right) }.%
\mathbf{P}_{k}^{\left( s+1\right) }\mathbf{P}_{l}^{\left( s+1\right)
}\right) _{\mathrm{M}}=\left\vert \mathbf{P}_{k}^{\left( s\right) }\mathbf{P}%
_{l}^{\left( s\right) }\right\vert _{\mathrm{M}}\left\vert \mathbf{P}%
_{k}^{\left( s+1\right) }\mathbf{P}_{l}^{\left( s+1\right) }\right\vert _{%
\mathrm{M}}\cosh \left( \Delta \varphi _{kl}^{\left( s\right) }\right)
\label{b5.10}
\end{equation}%
where $\Delta \varphi _{kl}^{\left( s\right) }$ is the angle between
timelike vectors $\mathbf{P}_{k}^{\left( s\right) }\mathbf{P}_{l}^{\left(
s\right) }$ and $\mathbf{P}_{k}^{\left( s+1\right) }\mathbf{P}_{l}^{\left(
s+1\right) }$ in the space-time geometry of Minkowski. The angle $\Delta
\varphi _{kl}^{\left( s\right) }$ may describe a rotation in the following
sense. If two vectors $\mathbf{P}_{0}\mathbf{P}_{1}$ and $\mathbf{P}_{0}%
\mathbf{P}_{2}$ have the common origin, and their length are equal $%
\left\vert \mathbf{P}_{0}\mathbf{P}_{1}\right\vert _{\mathrm{M}}=\left\vert
\mathbf{P}_{0}\mathbf{P}_{2}\right\vert $, the scalar product of two vectors
$\left( \mathbf{P}_{0}\mathbf{P}_{1}.\mathbf{P}_{0}\mathbf{P}_{1}\right) _{%
\mathrm{M}}$ can be presented in the form$_{\mathrm{M}}$%
\begin{equation}
\left( \mathbf{P}_{0}\mathbf{P}_{1}.\mathbf{P}_{0}\mathbf{P}_{2}\right) _{%
\mathrm{M}}=\left\vert \mathbf{P}_{0}\mathbf{P}_{1}\right\vert _{\mathrm{M}%
}\cdot \left\vert \mathbf{P}_{0}\mathbf{P}_{2}\right\vert _{\mathrm{M}}\cosh
\varphi  \label{b5.11}
\end{equation}%
where $\varphi $ is the angle between the vectors $\mathbf{P}_{0}\mathbf{P}%
_{1}$ and $\mathbf{P}_{0}\mathbf{P}_{2}$. Let $\varphi =0$, then%
\begin{equation}
\left( \mathbf{P}_{0}\mathbf{P}_{1}.\mathbf{P}_{0}\mathbf{P}_{2}\right) _{%
\mathrm{M}}=\left\vert \mathbf{P}_{0}\mathbf{P}_{1}\right\vert _{\mathrm{M}%
}\cdot \left\vert \mathbf{P}_{0}\mathbf{P}_{2}\right\vert _{\mathrm{M}}
\label{b5.12}
\end{equation}%
and vectors $\mathbf{P}_{0}\mathbf{P}_{1}$ and $\mathbf{P}_{0}\mathbf{P}_{2}$
coincide, i.e. $P_{1}=P_{2}$, if vectors $\mathbf{P}_{0}\mathbf{P}_{1}$ and $%
\mathbf{P}_{0}\mathbf{P}_{2}$ are timelike $\left\vert \mathbf{P}_{0}\mathbf{%
P}_{1}\right\vert _{\mathrm{M}}=\left\vert \mathbf{P}_{0}\mathbf{P}%
_{2}\right\vert _{\mathrm{M}}>0$. If vectors $\mathbf{P}_{0}\mathbf{P}_{1}$
and $\mathbf{P}_{0}\mathbf{P}_{2}$ are spacelike ($\left\vert \mathbf{P}_{0}%
\mathbf{P}_{1}\right\vert _{\mathrm{M}}=\left\vert \mathbf{P}_{0}\mathbf{P}%
_{2}\right\vert _{\mathrm{M}}<0$), then the vectors $\mathbf{P}_{0}\mathbf{P}%
_{1}$ and $\mathbf{P}_{0}\mathbf{P}_{2}$ do not coincide ($P_{1}\neq P_{2}$%
), in general, although the angle $\varphi $ between them vanishes. Such a
situation takes place, because the geometry of Minkowski is single-variant
with respect to timelike vectors, and it is multivariant with respect to
space-like vectors. The presentation of the scalar product in the form (\ref%
{b5.10}) is useful only for timelike vectors ($\left\vert \mathbf{P}%
_{k}^{\left( s\right) }\mathbf{P}_{l}^{\left( s\right) }\right\vert _{%
\mathrm{M}}>0,\left\vert \mathbf{P}_{k}^{\left( s+1\right) }\mathbf{P}%
_{l}^{\left( s+1\right) }\right\vert _{\mathrm{M}}>0$), because for
spacelike vectors the angle $\Delta \varphi _{kl}^{\left( s\right) }$ does
not describe completely the mutual disposition of vectors.

Let both vectors $\mathbf{P}_{k}^{\left( s\right) }\mathbf{P}_{l}^{\left(
s\right) }$ and $\mathbf{P}_{k}^{\left( s+1\right) }\mathbf{P}_{l}^{\left(
s+1\right) }$ are timelike in geometry of Minkowski, i.e.
\begin{equation}
\left\vert \mathbf{P}_{k}^{\left( s\right) }\mathbf{P}_{l}^{\left( s\right)
}\right\vert _{\mathrm{M}}>0,\qquad \left\vert \mathbf{P}_{k}^{\left(
s+1\right) }\mathbf{P}_{l}^{\left( s+1\right) }\right\vert _{\mathrm{M}}>0
\label{b5.10a}
\end{equation}%
By means of (\ref{b5.10}) the equations (\ref{b5.7}) can be rewritten in the
form
\begin{eqnarray}
&&\left\vert \mathbf{P}_{k}^{\left( s\right) }\mathbf{P}_{l}^{\left(
s\right) }\right\vert _{\mathrm{M}}\left\vert \mathbf{P}_{k}^{\left(
s+1\right) }\mathbf{P}_{l}^{\left( s+1\right) }\right\vert _{\mathrm{M}%
}\cosh \left( \Delta \varphi _{kl}^{\left( s\right) }\right) -2\sigma _{%
\mathrm{M}}\left( P_{k}^{\left( s\right) },P_{l}^{\left( s\right) }\right)
\notag \\
&=&2d\left( P_{k}^{\left( s\right) },P_{l}^{\left( s\right) }\right)
-w\left( P_{k}^{\left( s\right) },P_{l}^{\left( s\right) },P_{k}^{\left(
s+1\right) },P_{l}^{\left( s+1\right) }\right)   \label{b5.14}
\end{eqnarray}%
\begin{equation}
\cosh \left( \Delta \varphi _{kl}^{\left( s\right) }\right) =\frac{%
\left\vert \mathbf{P}_{k}^{\left( s\right) }\mathbf{P}_{l}^{\left( s\right)
}\right\vert _{\mathrm{M}}}{\left\vert \mathbf{P}_{k}^{\left( s+1\right) }%
\mathbf{P}_{l}^{\left( s+1\right) }\right\vert _{\mathrm{M}}}+\frac{2d\left(
P_{k}^{\left( s\right) },P_{l}^{\left( s\right) }\right) -w\left(
P_{k}^{\left( s\right) },P_{l}^{\left( s\right) },P_{k}^{\left( s+1\right)
},P_{l}^{\left( s+1\right) }\right) }{\left\vert \mathbf{P}_{k}^{\left(
s\right) }\mathbf{P}_{l}^{\left( s\right) }\right\vert _{\mathrm{M}%
}\left\vert \mathbf{P}_{k}^{\left( s+1\right) }\mathbf{P}_{l}^{\left(
s+1\right) }\right\vert _{\mathrm{M}}}  \label{b5.15}
\end{equation}%
In the case, when the space-time geometry is uniform and isotropic in the
sense that%
\begin{equation}
d\left( P_{k}^{\left( s\right) },P_{l}^{\left( s\right) }\right) =d\left(
\sigma _{\mathrm{M}}\left( P_{k}^{\left( s\right) },P_{l}^{\left( s\right)
}\right) \right) ,\qquad \forall P_{k}^{\left( s\right) },P_{l}^{\left(
s\right) }\in \Omega   \label{b5.16}
\end{equation}%
the following relations take place
\begin{equation}
\left\vert \mathbf{P}_{k}^{\left( s\right) }\mathbf{P}_{l}^{\left( s\right)
}\right\vert _{\mathrm{M}}=\left\vert \mathbf{P}_{k}^{\left( s+1\right) }%
\mathbf{P}_{l}^{\left( s+1\right) }\right\vert _{\mathrm{M}},\qquad
k,l=0,1,...n  \label{b5.17}
\end{equation}%
The dynamic equation (\ref{b5.15}) can be rewritten in the form%
\begin{equation}
\cosh \left( \Delta \varphi _{kl}^{\left( s\right) }\right) -1=2\sinh ^{2}%
\frac{\Delta \varphi _{kl}^{\left( s\right) }}{2}=\frac{2d\left(
P_{k}^{\left( s\right) },P_{l}^{\left( s\right) }\right) -w\left(
P_{k}^{\left( s\right) },P_{l}^{\left( s\right) },P_{k}^{\left( s+1\right)
},P_{l}^{\left( s+1\right) }\right) }{\left\vert \mathbf{P}_{k}^{\left(
s\right) }\mathbf{P}_{l}^{\left( s\right) }\right\vert _{\mathrm{M}}^{2}}
\label{b5.18}
\end{equation}%
where $w\left( P_{k}^{\left( s\right) },P_{l}^{\left( s\right)
},P_{k}^{\left( s+1\right) },P_{l}^{\left( s+1\right) }\right) $ is
determined by the relation (\ref{b5.9}).

In the case of the space-time geometry (\ref{b4.2}), where $d\left(
P_{0},P_{1}\right) =$const for timelike vectors $\mathbf{P}_{0}\mathbf{P}%
_{1} $, the dynamic equation (\ref{b5.18}) takes the form%
\begin{equation}
\sinh \frac{\Delta \varphi _{kl}^{\left( s\right) }}{2}=\sqrt{\frac{d}{%
2\sigma _{\mathrm{M}}\left( P_{k}^{\left( s\right) },P_{l}^{\left( s\right)
}\right) }},\qquad \sigma _{\mathrm{M}}\left( P_{k}^{\left( s\right)
},P_{l}^{\left( s\right) }\right) >0,\qquad k\neq l  \label{b5.19}
\end{equation}

In particular case of a pointlike particle, when the link of the world chain
consists of two points $\mathcal{P}^{1}=\left\{ P_{0},P_{1}\right\} $ the
dynamic equation (\ref{b5.19}) takes the form
\begin{equation}
\sinh \frac{\Delta \varphi _{01}^{\left( s\right) }}{2}=\sqrt{\frac{d}{%
2\sigma _{\mathrm{M}}\left( P_{0}^{\left( s\right) },P_{1}^{\left( s\right)
}\right) }}  \label{b5.20}
\end{equation}%
The world chain of the pointlike particle in the space-time geometry (\ref%
{b4.2}) with $d=\hbar /\left( 2bc\right) $ is described in terms of the
geometry of Minkowski. Any next link $\mathbf{P}_{0}^{\left( s+1\right) }%
\mathbf{P}_{1}^{\left( s+1\right) }$ of the world chain is obtained from the
previous link $\mathbf{P}_{0}^{\left( s\right) }\mathbf{P}_{1}^{\left(
s\right) }$ as a result of translation at the distance $\sqrt{2\sigma _{%
\mathrm{M}}\left( P_{0}^{\left( s\right) },P_{1}^{\left( s\right) }\right) }$
and rotation with the angle $\Delta \varphi _{kl}^{\left( s\right) }$, which
is defined by the dynamic equation (\ref{b5.20}). The value of the angle $%
\Delta \varphi _{kl}^{\left( s\right) }$ is determined, but the axis of
rotation remains to be indefinite. As a result position of the next link $%
\mathbf{P}_{0}^{\left( s+1\right) }\mathbf{P}_{1}^{\left( s+1\right) }$
appears to be indefinite. According to dynamic equation (\ref{b5.20})
possible positions of the link $\mathbf{P}_{0}^{\left( s+1\right) }\mathbf{P}%
_{1}^{\left( s+1\right) }$ form a cone with the angle $\Delta \varphi
_{kl}^{\left( s\right) }$ at it vertex. As far as such a cone appears at any
link of the world chain, the chain appears to be multivariant. Statistical
description of such wobbling chains leads to the Schr\"{o}dinger equation
\cite{R91} for a free particle of the mass
\begin{equation}
m=b^{-1}\sqrt{2\sigma \left( P_{0}^{\left( s\right) },P_{1}^{\left( s\right)
}\right) }=b^{-1}\sqrt{2\sigma _{\mathrm{M}}\left( P_{0}^{\left( s\right)
},P_{1}^{\left( s\right) }\right) +\frac{\hbar }{bc}}  \label{b5.21}
\end{equation}

In terms of the physical space-time geometry $\mathcal{G}$, described by the
world function $\sigma $, the world chain appears to be also multivariant,
although in terms of the world function $\sigma $ there is only translation
of links in the timelike direction without rotation. Nevertheless, the
translation appears to be multivariant, because the world function $\sigma $
describes space-time geometry, which is multivariant with respect to both
timelike and spacelike vectors.

Thus, dynamics of free composite particles may be described in any physical
geometry. Dynamics is multivariant, in general. It is determined completely
by the space-time geometry and does not refer to any means of description.
The free motion of the composite particle may be reduced to description in
the framework of the space-time geometry of Minkowski. However, in this case
the motion of a particle ceases to be free. One obtains the motion in the
field of some metric fields, generated by the space-time geometry. This
effect is well known in the general relativity. The free particle, moving in
the curved space-time, may be described as a particle, moving in the
gravitational field inside the space-time of Minkowski.

Dynamics of free composite particles is determined completely by the
space-time geometry and only by the space-time geometry. The mass of the
pointlike particle is geometrized in the sense, that it is determined by the
length of links of the world chain. Dynamic equations, represented in terms
of the geometry of Minkowski appear to be equations in finite differences
(but not differential equations). Under some condition they can be reduced
to differential equations. Dynamics is not restricted by the demands of the
space-time continuity and those of unlimited divisibility of the space-time.
The dynamics does not contain a reference to the linear vector space and its
properties, although the dynamics may be described in terms of the
space-time geometry of Minkowski, which uses properties of the linear vector
space.

\section{Concluding remarks}

The geometry with the intransitive equivalence relation is a new conception
of geometry. Its name has not been established yet. There are different
names of this geometry, reflecting the process of its development and its
various properties: physical geometry, tubular geometry (T-geometry),
multivariant geometry, intransitive geometry. The multivariant geometry is
constructed on the basis of the deformation principle, when all propositions
of the multivariant geometry are obtained as a result of deformation of the
proper Euclidean geometry. This method of construction is free of defects of
the conventional Euclidean method, when all geometric proposition are
deduced from some axiomatics. The Euclidean method is not applicable for
construction of a physical geometry, because almost all physical geometries
are not axiomatizable. For instance, the Riemannian geometry cannot be
axiomatized. If nevertheless the Riemannian geometry is constructed by the
conventional method on the basis of axioms, which are not formulated
explicitly, it appears to be inconsistent. It differs in some details from
the $\sigma $-Riemannian geometry, i.e. from the Riemannian geometry,
constructed by means of the deformation principle.

The application of the deformation principle for construction of a physical
geometry is free of such problems as a test of the axiomatics consistency
and deduction of geometric propositions, accompanied by formulation and
proof of numerous theorems.

The class of physical geometries is more powerful, then the class of the
axiomatizable geometries, which are logical constructions. These logical
constructions (excluding the proper Euclidean geometry) have not a direct
relation to the geometry, as a science on mutual disposition of geometrical
objects in the space and in the space-time.

The physical geometry admits one to consider and describe such properties of
the space-time as discreteness, granularity and limited divisibility of the
space-time. These properties cannot be described by the geometries,
constructed on the basis on the linear vector space and unlimited
divisibility of the space. The intransitive equivalence relation and the
connected with it multivariance are a source of the granularity and limited
divisibility of the space-time. These properties of the space-time geometry
are important for construction of the discrimination mechanism, generating
discrete parameters of elementary particles. The physical space-time
geometry admits one to explain the elementary particles structure
geometrically, i.e. as a simple combinations of space-time points, ignoring
such secondary constructions as wave functions and other quantum attributes.

The physical geometry is a further development of geometrical methods in its
application to investigation of the space-time properties. Application of
the physical geometry methods to investigation of the space-time realizes
further geometrizaition of the microcosm physics, which was made start by
the spacial relativity and the general relativity.

\end{document}